\documentclass[12pt]{article}
\usepackage{graphicx} %% Package for Figure
\usepackage{float} %% Package for Float
\usepackage{amssymb}
\usepackage{biblatex}
\usepackage{amsmath}
\usepackage{indentfirst}
\usepackage{mathtools}
\numberwithin{equation}{section}

\newtheorem{lemma}{Lemma}[section]

\graphicspath{{figure/}}                                 %% Path of Figures
\usepackage[a4paper]{geometry}                           %% Paper size
\begin{document}
\title{\bf Boundedness of Forelli-Rudin Type Operators on Tubular Domains over The Generalized Light Cones
\thanks{The work was Supported by NSFC (11971042, 12071035) , the National Key R\&D Program of China (2021YFA1002600) and Natural Science Foundation of Beijing Municipal (No.1252005).\vskip 1mm
Corresponding author's E-mail address: denggt@bnu.edu.cn }}
\author{\small Xin Xia $^{1}$. GuanTie Deng$^{2}$\\
 \small 1.2\,\,School of Mathematical Sciences, Beijing Normal University, Beijing 100875, China
}
\date{}
\maketitle
\noindent{\bf Abstract:} This study investigates conditions for the boundedness of Forelli-Rudin type operators on weighted Lebesgue spaces associated with tubular domains over the generalized light cone. We establish a complete characterization of the boundedness for two classes of Forelli-Rudin type operators from  $L_{\boldsymbol{\alpha}}^{p}$ to $L_{\boldsymbol{\beta}}^{q}$, in the range $1 <  p \leq q < \infty$. The findings contribute significantly to the analysis of Bergman projection operators in this setting.   

\section{Introduction and main results}
In this article, we study two classes of Forelli-Rudin type operators acting on weighted Lebesgue spaces. Our main objective is to establish boundedness conditions for such operators in the setting of tubular domains over the generalized light cone. To this end, we first introduce relevant definitions and notations.

Let $\mathbb{C}^n$ denote the $n$-dimensional complex Euclidean space, where $n$ is a positive integer.  For any two points \(z=(z_{1},\ldots,z_{n})\) and \(w=(w_{1},\ldots,w_{n})\) in \(\mathbb{C}^{n}\), we write
\[
\langle z, w\rangle:=z_{1}\overline{w}_{1}+\cdots+z_{n}\overline{w}_{n}
\]
and \(|z|:=\sqrt{\langle z, z\rangle}\).

For \(z\in\mathbb{C}^{n}\), we also use the notation
\[
z=(z',z_{n}),\quad\text{where }z'=(z_{1},\ldots,z_{n - 1})\in\mathbb{C}^{n - 1}\text{ and }z_{n}\in\mathbb{C}^{1}.
\]

Let \(\mathbb{B}_{n}:=\left\{z\in \mathbb{C}^{n}:\vert z\vert< 1\right\}\) be the open unit ball in \(\mathbb{C}^{n}\) and \(S_{n}:=\left\{z\in \mathbb{C}^{n}:\vert z\vert = 1\right\}\) denote its boundary. For \(a,b,c\in \mathbb{R}\), the Forelli - Rudin type operators are defined as

\[T_{a,b,c}f(z)=(1-\vert z\vert^{2})^{a}\int_{\mathbb{B}_{n}}\frac{(1-\vert u\vert^{2})^{b}f(u)d\nu(u)}{(1-\langle z, u\rangle)^{c}}\]

and

\[S_{a,b,c}f(z)=(1-\vert z\vert^{2})^{a}\int_{\mathbb{B}_{n}}\frac{(1-\vert u\vert^{2})^{b}f(u)d\nu(u)}{\vert1-\langle z, u\rangle\vert^{c}},\]
where $d\nu$ is the volume measure on $\mathbb{B}_{n}$, normalized so that $\nu(\mathbb{B}_{n}) = 1$. Also, for any real parameter $\alpha$ we define $d\nu_{\alpha}(z):=(1 - |z|^{2})^{\alpha}d\nu(z)$. The study of these operators in special cases dates back to Stein \([15]\), who established the boundedness of \(T_{0,0,n+1}\) on \(L^p(\mathbb{B}_n)\) for \(1 < p < \infty\). Kolaski \([7]\) later revisited this operator from the perspective of Bergman projections. In 1991, Zhu \([18]\) derived the necessary and sufficient condition for the boundedness of the general operator \(T_{a,b,c}\) in the one-dimensional case (\(n = 1\)) with \(c = 1 + a + b\).
This result has since been extended to higher dimensions by Kures and Zhu \([6]\), who established the following two theorems.

\textbf{Theorem A} Suppose $1 < p < \infty$. Then the following conditions are equivalent:
\begin{enumerate}
    \item[(i)] The operator $T_{a,b,c}$ is bounded on $L^{p}(\mathbb{B}_n,d\nu_{\alpha})$.
     \item[(ii)] The operator $S_{a,b,c}$ is bounded on $L^{p}(\mathbb{B}_n,d\nu_{\alpha})$.
     \item[(iii)] The parameters satisfy
    \[
    \begin{cases}
    -pa < \alpha + 1 < p(b + 1)\\
    c\leq n + 1+a + b.
    \end{cases}
    \]
\end{enumerate}

\textbf{Theorem B} The following conditions are equivalent:
\begin{enumerate}
    \item[(i)] The operator $T_{a,b,c}$ is bounded on $L^{1}(\mathbb{B}_n,d\nu_{\alpha})$.
    \item[(ii)] The operator $S_{a,b,c}$ is bounded on $L^{1}(\mathbb{B}_n,d\nu_{\alpha})$.
    \item[(iii)] The parameters satisfy
    \[
    \begin{cases}
    -a < \alpha + 1 < b + 1\\
    c=n + 1+a + b
    \end{cases}
    \text{ or }
    \begin{cases}
    -a < \alpha + 1\leq b + 1\\
    c< n + 1+a + b.
    \end{cases}
    \]
\end{enumerate}

These two theorems were originally established in [6] under the additional assumption that \( c \) is neither zero nor a negative integer. Recently, Zhao [23] removed this restriction and further extended the results by characterizing the boundedness of both \( T_{a,b,c} \) and \( S_{a,b,c} \) from \( L^p(\mathbb{B}, dv_{\alpha}) \) to \( L^q(\mathbb{B}, dv_{\beta}) \), where \( 1 \leq p \leq q < \infty \).

For our purpose, we introduce some definitions and notations. Let \(\Omega\) be an arbitrary subset of \(\mathbb{R}^{n}\). The tubular region \(T_{\Omega}\) is defined as
\[T_{\Omega}=\left\{z=x + iy\in\mathbb{C}^{n}:x\in\mathbb{R}^{n},y\in\Omega\right\}.\]

 On the tubular region \(T_{\Omega}\), we define the function \(\rho(z)=\rho(iy)\) to be a positive continuous function and use the function \(\rho(z)\) to define the weighted function space. Define the weighted Lebesgue space \(L_{\alpha}^{p}(T_{\Omega})\): assume the measure
\[dV_{\alpha}(z)=\rho(z)^{\alpha}dV(z)\]
where \(dV = dxdy\) is the Lebesgue measure, and \(\alpha\) is a real parameter. For \(1\leqslant p<\infty\) and \(\alpha\in\mathbb{R}\), the weighted Lebesgue norm $\|f\|_{L_{\alpha}^{p}}$ is
$$\|f\|_{L_{\alpha}^{p}}=\left(\int_{T_{\Omega}}|f(z)|^{p}dV_{\rho}(z)\right)^{\frac{1}{p}}.$$

we suppose
$$X=\left(
\begin{array}{cccc}
x_1 & 0 & \cdots & x_{2n - 1}\\
0 & x_{2} & 0 & x_{2n-2}\\
\vdots & 0 & \ddots & \vdots\\
x_{2n - 1} & x_{2n-2} & x_{n+1} & x_{n}
\end{array}
\right),$$
the generalized light cone \(\mathcal{P}_{n}\) is defined as
$$\mathcal{P}_n = \{x \in \mathbb{R}^m : x_{1}>0,x_{2}>0,\cdots,x_{n-1}>0,det(x)>0\}.$$

In addition, we can define the function $\Delta_{k}(Y)$ as the $k$-th leading principal minor of the matrix $Y$, where $k = 1, 2,\cdots, n$, that is
\[
\Delta_{k}(Y)=\begin{vmatrix}
y_{11}&\cdots&y_{1k}\\
\vdots&\ddots&\vdots\\
y_{1k}&\cdots&y_{kk}
\end{vmatrix}.
\]

If $x\in\mathcal{P}_{n}$, then $\Delta_{k}(X)=x_{1}x_{2}\cdots x_{k}$, $k = 1, 2,\cdots, n-1$ and $$\Delta_{n}(X)=det(X)=x_{1}x_{2}\cdots x_{n-1}(x_{n}-\frac{x_{2n-1}^{2}}{x_{1}}-\cdots-\frac{x_{n+1}^{2}}{x_{n-1}})$$
by definition, we know  \(X\) is a positive-definite matrix.

That is, \(\mathcal{P}_{n}\) is a convex open set in \(\mathbb{R}^{m}\) $(m=2n-1)$, we may regard the elements in \(\mathcal{P}_{n}\) as matrices and perform operations with the same operation rules as matrix operations.

For \( y \in \mathcal{P}_n \) and \( \boldsymbol{\alpha} = (\alpha_1, \alpha_2, \cdots, \alpha_n) \), Define a function
\[
\Delta^{\boldsymbol{\alpha}}(y) = \Delta_1(y)^{\alpha_1 - \alpha_2} \Delta_2(y)^{\alpha_2 - \alpha_3} \cdots \Delta_{n-1}(y)^{\alpha_{n-1} - \alpha_n} \Delta_n(y)^{\alpha_n}.
\]

The functions \( P^{\boldsymbol{\alpha}}(z) \) is defined as
\(
P^{\boldsymbol{\alpha}}(z) = \Delta^{\boldsymbol{\alpha}}\left( \frac{x + iy}{i} \right)
\), the measure \( dV_{\boldsymbol{\alpha}}(z) = \Delta^{\boldsymbol{\alpha}}(y) dV(z) \), where
\[
dV(z) = dxdy = \prod_{i \leqslant j} dx_{ij} \cdot \prod_{i \leqslant j} dy_{ij}.
\]

For $f\in$ \(L^{p}_{\alpha}(T_{\mathcal{P}_{n}})\), the operator \(T\) is defined in this chapter as
\[Tf(z):=T^{\mathcal{P}_{n}}_{\mathbf{a}, \mathbf{b}, \mathbf{c}}f(z)=\Delta^{\mathbf{a}}(\text{Im}z)\int_{T_{\Lambda_{n}}}\frac{\Delta^{\mathbf{b}}(\text{Im}w)}{P^{\mathbf{c}}(z-\overline{w})}f(w)dV(w)\]
where \begin{align*}\begin{cases}\mathbf{a}=(a_1+\frac{n-2}{2},a_2
+\frac{n-2}{2},\cdots,a_n),\\
\mathbf{b}=(b_1+\frac{n-2}{2},b_2+\frac{n-2}{2},\cdots,b_n),\\
\mathbf{c}=(c_1+\frac{n-2}{2},c_2+\frac{n-2}{2},\cdots,c_n).
\end{cases}
\end{align*}

Our main results are stated as:

\textbf{Theorem 1} Suppose that \( 1 < p \leq q < \infty \). If the operator \( T \) is a bounded mapping from \(L_{\boldsymbol{\alpha}}^{p}(\mathcal{P}_{n}) \) to \( L_{\boldsymbol{\beta}}^{q}(\mathcal{P}_{n}) \), then the parameters \( \mathbf{a} \), \( \mathbf{b} \), and \( \mathbf{c} \) should satisfy for all \(j = 1,2,\cdots,n-1\)
\begin{align*}\begin{cases}
-a_jq< \beta_j+\frac{n+1}{2},\quad \alpha_j+\frac{n+1}{2} < p(b_j+\frac{n+1}{2}),\\
c_j = a_j + b_j + n + 1 + \frac{\beta_j + n + 1}{q} - \frac{\alpha_j + n + 1}{p},
\end{cases}
\end{align*}
and
\begin{align*}\begin{cases}
-a_nq< \beta_n+1,\quad\alpha_n+1 < p\left(b_n+1\right),\\
c_n = a_n + b_n + n + 1 + \frac{\beta_n + n + 1}{q} - \frac{\alpha_n + n + 1}{p}.
\end{cases}
\end{align*}

\textbf{Theorem 2} Suppose that \( 1 < p \leq q < \infty \). If for all \(j = 1,2,\cdots,n-1\), the parameters satisfy the condition \( c_j > n \),
\begin{align*}\begin{cases}
\alpha_j + 1 < p\left( \left( \frac{1}{2q} + \frac{1}{2p} \right)(1 - n) + b_j + \frac{n + 1}{2} \right), \\
\beta_j + 1 > q\left( \left( \frac{1}{2q} + \frac{1}{2p} \right)(1 - n) - a_j \right),\\
c_j = a_j + b_j + n + 1 + \frac{\beta_j + n + 1}{q} - \frac{\alpha_j + n + 1}{p},
\end{cases}
\end{align*}
 and \( c_n > n \),
 \begin{align*}\begin{cases}
 \alpha_n + 1 < p( b_n + 1),\\
 -a_nq< \beta_n+1,\\
 c_n = a_n + b_n + n + 1 + \frac{\beta_n + n + 1}{q} - \frac{\alpha_n + n + 1}{p}.
 \end{cases}
\end{align*}

Then the operator \( T \) is bounded from \(L_{\boldsymbol{\alpha}}^{p}(\mathcal{P}_{n}) \)  to $L_{\boldsymbol{\beta}}^{q}(\mathcal{P}_{n}) $.

\section{Preliminaries}
\begin{lemma}$[2]$ Let \(\Omega\) be a domain in \(\mathbb{R}^n\). Then the reproducing kernel in the weighted Bergman space \(A_{\alpha}^{p}(T_{\Omega})\) can be expressed as:
\[
K(z, w)=\int_{\mathbb{R}^n}e^{2\pi i(z - \overline{w})\cdot t}I^{-1}(t)dt,
\]
where \(I(t)=\int_{\Omega}\rho(iy)^{\alpha}e^{-4\pi y\cdot t}dy\).

\end{lemma}

\begin{lemma} $[22]$ Let \(p,r,q\) be positive numbers satisfying \(1<p\leqslant r\) and \(\frac{1}{p}+\frac{1}{q}=1\)
Let \(H(x,y)\) be a non - negative measurable function on \(X\times Y\). Suppose there exist \(0 < t\leqslant1\), measurable functions \(\varphi_{1}:X\to(0,\infty)\), \(\varphi_{2}:Y\to(0,\infty)\) and non - negative constants \(M_{1},M_{2}\) such that
\[\int_{X}\int_{Y}H(x,y)^{t}\varphi_{1}^{q}(y)d\mu(y)\leqslant M_{1}^{q}\varphi_{2}^{q}(x)\]
almost everywhere on \(Y\) and
\[\int_{Y}H(x,y)^{(1 - t)r}\varphi_{2}^{r}(x)d\mu(x)\leqslant M_{2}^{r}\varphi_{1}^{r}(y)\]
almost everywhere on \(X\)

If
\[Tf(x)=\int_{X}f(y)H(x,y)d\mu(y)\]
where \(f\in L^{p}(X,d\mu)\), then \(T:L^{p}(X,d\mu)\to L^{r}(Y,d\nu)\) is bounded, and for each \(f\in L^{p}(X,d\mu)\)
\[\|Tf\|_{L^{r}(Y,d\nu)}\leqslant M_{1}M_{2}\|f\|_{L^{p}(X,d\mu)}.\]
\end{lemma}

\begin{lemma}
Suppose $
\boldsymbol{s}=(s_1,s_2,\cdots,s_n)\in\mathbb{R}^n, t\in\mathbb{R}^m
$, then we have the following properties:

(1) $\text{ If } t\in\mathcal{P}_{n}\text{ and }s_n>-1,s_{j}>-\frac{3}{2}, j = 1,2,\cdots,n-1$, the integral \begin{align*}
I_{\boldsymbol{s}}(t)&=\int_{\mathcal{P}_{n}}e^{- 4\pi y\cdot t}\Delta^{\boldsymbol{s}}(y)dy\\
&=C_{1,\boldsymbol{s}}t_{n}^{-s_{n}-\frac{n+1}{2}} \prod_{j=1}^{n-1}(4t_{j}-\frac{t_{2n-j}^{2}}{t_{n}})^{-s_{j}-\frac{3}{2}},
\end{align*}
where $$C_{1,\boldsymbol{s}}=\frac{\Gamma(s_{n}+1)\prod_{j=1}^{n-1}\Gamma(s_{j}+\frac{3}{2})}{2^{2s_{n}+n+1}\pi^{\sum_{j=1}^{n}s_{j}+\frac{3n-1}{2}}} .$$

\noindent \textbf{Proof}: For $t\in\mathcal{P}_{n}$ and $s_n>-1,s_{j}>-\frac{3}{2}\quad j = 1,2,\cdots,n-1$,
.
\begin{align*}
I_{\boldsymbol{s}}(t)=&\int_{\mathcal{P}_{n}}e^{- 4\pi y\cdot t}\Delta^{\boldsymbol{s}}(y)dy\\
=&\int_{\mathcal{P}_{n}}e^{-4\pi y\cdot t}\Delta_1(y)^{s_1 - s_2}\Delta_2(y)^{s_2 - s_3}\cdots\Delta_{n - 1}(y)^{s_{n - 1} - s_n}\Delta_n(y)^{s_n}dy\\
=&\int_{y_{1}>0}e^{-4\pi y_{1}t_{1}}(y_{1})^{s_{1}}dy_{1}\cdots\int_{y_{n-1}>0}e^{-4\pi y_{n-1}t_{n-1}}
(y_{n-1})^{s_{n-1}}dy_{n-1}\\
&\cdot\int_{\mathbb{R}^{n-1}}\int_{y_{n}-\frac{y_{2n-1}^{2}}{y_{1}}-\cdots-\frac{y_{n+1}^{2}}{y_{n-1}}>0}
e^{-4\pi (y_{n}t_{n}+y_{n+1}t_{n+1}+\cdots+y_{2n-1}t_{2n-1})}\\
&(y_{n}-\frac{y_{2n-1}^{2}}{y_{1}}-\cdots-\frac{y_{n+1}^{2}}{y_{n-1}})^{s_{n}}dy_{n}dy_{n+1}\cdots dy_{2n-1}.
\end{align*}

Let:
$u_i = y_i, v_j = y_{2n-j}$ for $j = 1, \dots, n-1$,$w = y_n$, the integral becomes:
\[
\begin{aligned}
I_{\boldsymbol{s}}(t) = & \left( \prod_{i=1}^{n-1} \int_0^\infty e^{-4\pi u_i t_i} u_i^{s_i}  du_i \right) \\
& \cdot \int_{\mathbb{R}^{n-1}} \int_{w > \sum_{j=1}^{n-1} \frac{v_j^2}{u_j}} e^{-4\pi (w t_n + \sum_{j=1}^{n-1} v_j t_{2n-j})} \left( w - \sum_{j=1}^{n-1} \frac{v_j^2}{u_j} \right)^{s_n}  dw  dv_1 \cdots dv_{n-1}.
\end{aligned}
\]

We set:
\[
z = w - \sum_{j=1}^{n-1} \frac{v_j^2}{u_j} \quad \Rightarrow \quad w = z + \sum_{j=1}^{n-1} \frac{v_j^2}{u_j}, \quad dz = dw.
\]

The $w$-integral becomes:
\[
\int_0^\infty e^{-4\pi t_n z} z^{s_n}  dz \cdot e^{-4\pi t_n \sum_{j=1}^{n-1} \frac{v_j^2}{u_j}} \cdot e^{-4\pi \sum_{j=1}^{n-1} v_j t_{2n-j}}.
\]

We have:
\[
\int_0^\infty e^{-4\pi t_n z} z^{s_n}  dz = \frac{\Gamma(s_n+1)}{(4\pi t_n)^{s_n+1}}.
\]

So the inner integral over $w$ yields:
\[
\frac{\Gamma(s_n+1)}{(4\pi t_n)^{s_n+1}} \cdot \prod_{j=1}^{n-1} e^{-4\pi t_n \frac{v_j^2}{u_j} - 4\pi v_j t_{2n-j}}.
\]

For each $j$, compute:
\[
J_j = \int_{-\infty}^\infty e^{-4\pi t_n \frac{v_j^2}{u_j} - 4\pi v_j t_{2n-j}}  dv_j.
\]

Now:
\[
I_{\boldsymbol{s}}(t) = \frac{\Gamma(s_n+1)}{(4\pi t_n)^{s_n+1}} \cdot \prod_{j=1}^{n-1} \int_0^\infty e^{-4\pi u_j t_j} u_j^{s_j} \cdot J_j  du_j.
\]

Due to

\[
\begin{aligned}
-4\pi t_n \frac{v_j^2}{u_j} - 4\pi v_j t_{2n-j} &= -\frac{4\pi t_n}{u_j} \left( v_j^2 + \frac{u_j t_{2n-j}}{t_n} v_j \right) \\
&= -\frac{4\pi t_n}{u_j} \left[ \left( v_j + \frac{u_j t_{2n-j}}{2t_n} \right)^2 - \frac{u_j^2 t_{2n-j}^2}{4t_n^2} \right] \\
&= -\frac{4\pi t_n}{u_j} \left( v_j + \frac{u_j t_{2n-j}}{2t_n} \right)^2 + \frac{\pi u_j t_{2n-j}^2}{t_n}.
\end{aligned}
\]

Then
\[
\begin{aligned}
J_j &= e^{\frac{\pi u_j t_{2n-j}^2}{t_n}} \int_{-\infty}^\infty e^{-\frac{4\pi t_n}{u_j} v^2} dv\\
&= e^{\frac{\pi u_j t_{2n-j}^2}{t_n}} \cdot \sqrt{\frac{u_j}{4 t_n}}.
\end{aligned}
\]

That is:
\[
\begin{aligned}
I_j &= \int_0^\infty e^{-4\pi u_j t_j} u_j^{s_j} \cdot \sqrt{\frac{u_j}{4 t_n}} \exp\left( \frac{\pi t_{2n-j}^2 u_j}{t_n} \right)  du_j \\
&= \frac{1}{2\sqrt{t_n}} \int_0^\infty u_j^{s_j + \frac12} \exp\left[ -u_j \left( 4\pi t_j - \frac{\pi t_{2n-j}^2}{t_n} \right) \right] du_j.
\end{aligned}
\]

Let:
\[
\beta_j = 4\pi t_j - \frac{\pi t_{2n-j}^2}{t_n}.
\]

Then:
\[
I_j = \frac{1}{2\sqrt{t_n}} \cdot \frac{\Gamma(s_j + \frac32)}{\beta_j^{s_j + \frac32}}.
\]

Substitute $\beta_j$:
\[
\beta_j = \frac{\pi}{t_n} (4 t_n t_j - t_{2n-j}^2).
\]

So:
\[
\begin{aligned}
I_j &= \frac{1}{2\sqrt{t_n}} \cdot \frac{\Gamma(s_j + \frac32)}{\left( \frac{\pi}{t_n} (4 t_n t_j - t_{2n-j}^2) \right)^{s_j + \frac32}} \\
&= \frac{\Gamma(s_j + \frac32)}{2\pi^{s_j + \frac32}} \cdot \frac{t_n^{s_j + 1}}{(4 t_n t_j - t_{2n-j}^2)^{s_j + \frac32}}.
\end{aligned}
\]

We have:
\[
\begin{aligned}
I_{\boldsymbol{s}}(t) &= \frac{\Gamma(s_n+1)}{(4\pi t_n)^{s_n+1}} \cdot \prod_{j=1}^{n-1} I_j \\
&= \frac{\Gamma(s_n+1)}{(4\pi t_n)^{s_n+1}} \cdot \prod_{j=1}^{n-1} \left[ \frac{\Gamma(s_j + \frac32)}{2\pi^{s_j + \frac32}} \cdot \frac{t_n^{s_j + 1}}{(4 t_n t_j - t_{2n-j}^2)^{s_j + \frac32}} \right].
\end{aligned}
\]

Thus:
\[
\begin{aligned}
I_{\boldsymbol{s}}(t) = &\ \Gamma(s_n+1) \prod_{j=1}^{n-1} \Gamma(s_j + \tfrac32) \cdot \frac{1}{2^{2s_n + 2} \pi^{s_n+1} t_n^{s_n+1}} \cdot \frac{1}{2^{n-1} \pi^{\sum_{j=1}^{n-1} s_j + \frac{3(n-1)}{2}}} \\
& \cdot t_n^{\sum_{j=1}^{n-1} s_j + n - 1} \cdot \prod_{j=1}^{n-1} (4 t_n t_j - t_{2n-j}^2)^{-s_j - \frac32}.\\
=&\frac{\Gamma(s_{n}+1)\prod_{j=1}^{n-1}\Gamma(s_{j}+\frac{3}{2})}{2^{2s_{n}+n+1}\pi^{\sum_{j=1}^{n}s_{j}+\frac{3n-1}{2}}} t_{n}^{\sum_{j=1}^{n-1}s_{j}-s_{n}+n-2} \prod_{j=1}^{n-1}(4t_{n}t_{j}-t_{2n-j}^{2})^{-s_{j}-\frac{3}{2}}\\
=&C_{1,\boldsymbol{s}}t_{n}^{-s_{n}-\frac{n+1}{2}} \prod_{j=1}^{n-1}(4t_{j}-\frac{t_{2n-j}^{2}}{t_{n}})^{-s_{j}-\frac{3}{2}}.
\end{aligned}
\]

 Let $\Delta^{\boldsymbol{s}}(t_{\delta}^{-1})$ is the power function on the cone $\mathcal{P}_n$, defined as:
$$
\Delta^{\boldsymbol{s}}(t_{\delta}^{-1}) = t_n^{s_n} \prod_{j=1}^{n-1} \left(4t_j - \frac{t_{2n-j}^{2}}{t_n}\right)^{s_j}.
$$

Thus
$$I_{\boldsymbol{s}}(t)=C_{1,\boldsymbol{s}}\Delta^{\boldsymbol{s}}(t_{\delta}^{-1})
\Delta(t_{\delta})^{-\frac{n+1}{2}}\Delta_{n-1}(t_{\delta})^{\frac{n-2}{2}}.$$

(2) For \(z=x + iy\in T_{\mathcal{P}_{n}}\), $s_n>-n-1,s_{j}>-\frac{5}{2}, j = 1,2,\cdots,n-1$, then
\[
\int_{\mathcal{P}_{n}} e^{2 \pi i z \cdot t} {I_{\boldsymbol{s}}(t)}^{-1} dt = C_{2,\boldsymbol{s}} P^{\boldsymbol{-s}}(z) P(z)^{-n-1}P_{n-1}(z)^{n-2}
\]
where \[
C_{2,\boldsymbol{s}} = 2^{\sum_{j=1}^{n} s_j + n - 1} \pi^{-2n + 1} \frac{\Gamma(s_n + n + 1)}{\Gamma(s_n + 1)} \prod_{j=1}^{n-1} \frac{\Gamma\left(s_j + \frac{5}{2}\right)}{\Gamma\left(s_j + \frac{3}{2}\right)}.
\]

\noindent \textbf{Proof}: For \(z=x + iy\in T_{\mathcal{P}_{n}}\), let

\[
\begin{aligned}
K(z)=&\int_{\mathcal{P}_{n}} e^{2 \pi i z \cdot t} {I_{\boldsymbol{s}}(t)}^{-1} dt\\
=&C_{1,\boldsymbol{s}}^{-1}\int_{\mathcal{P}_{n}} e^{2 \pi i z \cdot t} t_{n}^{s_{n}+\frac{n+1}{2}} \prod_{j=1}^{n-1}(4t_{j}-\frac{t_{n+j}^{2}}{t_{n}})^{s_{j}+\frac{3}{2}} dt
\end{aligned}
\]

Then
\begin{align*}
K(iy)&=\int_{\mathcal{P}_{n}} e^{-2 \pi  y \cdot t} t_{n}^{s_{n}+\frac{n+1}{2}} \prod_{j=1}^{n-1}(4t_{j}-\frac{t_{2n-j}^{2}}{t_{n}})^{s_{j}+\frac{3}{2}} dt\\
&=\int_{\mathcal{P}_{n}} e^{-2 \pi  (y_{1}t_{1}+\cdots +y_{n}t_{n}+y_{n+1}t_{n+1}+\cdots+y_{2n-1}t_{2n-1})} t_{n}^{s_{n}+\frac{n+1}{2}} \prod_{j=1}^{n-1}(4t_{j}-\frac{t_{2n-j}^{2}}{t_{n}})^{s_{j}+\frac{3}{2}}dt.\\
&=\int_{t_{1}>0}e^{-2\pi y_{1}t_{1}}(4t_{1}-\frac{t_{2n-1}^{2}}{t_{n}})^{s_{1}+\frac{3}{2}}dt_{1}\int_{t_{2}>0}e^{-2\pi y_{2}t_{2}}(4t_{2}-\frac{t_{2n-2}^{2}}{t_{n}})^{s_{2}+\frac{3}{2}}dt_{2}\\
&\cdots\int_{\mathbb{R}^{n-1}}\int_{t_{n}-\frac{t_{2n-1}^{2}}{t_{1}}-\cdots-\frac{t_{n+1}^{2}}{t_{n-1}}>0}e^{-2\pi (y_{n}t_{n}+y_{n+1}t_{n+1}+\cdots+y_{2n-1}t_{2n-1})}
(t_{n})^{s_{n}+\frac{n+1}{2}}dt_{n}dt_{n+1}\cdots dt_{2n-1}\\
&=2^{\sum_{j=1}^{n} s_j + n - 1} \pi^{-2n + 1} \frac{\Gamma(s_n + n + 1)}{\Gamma(s_n + 1)} \prod_{j=1}^{n-1} \frac{\Gamma\left(s_j + \frac{5}{2}\right)}{\Gamma\left(s_j + \frac{3}{2}\right)}
\prod_{j=1}^{n-1} y_j^{-s_j - 3} \left( y_n - \sum_{j=1}^{n-1} \frac{y_{2n-j}^2}{y_j} \right)^{-s_n - n - 1}.
\end{align*}

Morera's theorem ensures that \( K(z) \) is a holomorphic function. For \( R > 0 \) and \( |x| < R \), performing a Taylor expansion on \( K(x + iy) \), we have
\[
K(x + iy) = \sum_{k_1 = 0}^{\infty} \cdots \sum_{k_n = 0}^{\infty} \frac{1}{k_1! \cdots k_n!} \left. \frac{\partial^{k_1 + \cdots + k_n} K(z)}{\partial z_1^{k_1} \cdots \partial z_n^{k_n}} \right|_{z = iy} x_1^{k_1} \cdots x_n^{k_n}.
\]

According to the uniqueness theorem of analytic functions, we have
\[
K(z) = C_{2,\boldsymbol{s}} P^{\boldsymbol{-s}}(z) P(z)^{-n-1}P_{n-1}(z)^{n-2} .
\]

 \hfill $\blacksquare$

\end{lemma}

\noindent \textbf{Corollary 1}: Suppose \( \mathbf{s}' = (s_1+\frac{n-2}{2}, s_2+\frac{n-2}{2}, \ldots, s_n) \in \mathbb{R}^n \).

(1)If $t\in\mathcal{P}_{n}\text{ and }s_n>-1,s_{j}>-\frac{n+1}{2}, j = 1,2,\cdots,n-1$, then the integral
 \begin{align*}
I_{\boldsymbol{s}'}(t)&=\int_{\mathcal{P}_{n}}e^{- 4\pi y\cdot t}\Delta^{\boldsymbol{s}'}(y)dy\\
&=\frac{\Gamma(s_{n}+1)\prod_{j=1}^{n-1}\Gamma(s_{j}+\frac{n+1}{2})}{2^{2s_{n}+n+1}\pi^{\sum_{j=1}^{n}s_{j}+\frac{n^{2}+1}{2}}} t_{n}^{\sum_{j=1}^{n-1}s_{j}-s_{n}+\frac{n^{2}-n-2}{2}} \prod_{j=1}^{n-1}(4t_{n}t_{j}-t_{2n-j}^{2})^{-s_{j}-\frac{n+1}{2}}\\
&=C_{3,\boldsymbol{s}'}t_{n}^{-s_{n}-\frac{n+1}{2}} \prod_{j=1}^{n-1}(4t_{j}-\frac{t_{2n-j}^{2}}{t_{n}})^{-s_{j}-\frac{n+1}{2}}.\\
&=C_{3,\boldsymbol{s}'}\Delta^{\boldsymbol{s}'}(t_{\delta}^{-1})
\Delta(t_{\delta})^{-\frac{n+1}{2}}.
\end{align*}

(2) For \(z=x + iy\in T_{\mathcal{P}_{n}}\), $s_n>-n-1,s_{j}>-\frac{n+3}{2},j = 1,2,\cdots,n-1$, then
\[
\int_{\mathcal{P}_{n}} e^{2 \pi i z \cdot t} {I_{\boldsymbol{s}'}(t)}^{-1} dt = C_{4,\boldsymbol{s}'} P^{\boldsymbol{-s}'}(z) P(z)^{-n-1}
\]
where \[
C_{4,\boldsymbol{s}} = 2^{\sum_{j=1}^{n} s_j + \frac{n(n - 1)}{2}} \pi^{-2n + 1} \frac{\Gamma(s_n + n + 1)}{\Gamma(s_n + 1)} \prod_{j=1}^{n-1} \frac{\Gamma\left(s_j + \frac{n+3}{2}\right)}{\Gamma\left(s_j + \frac{n+1}{2}\right)}.
\]

\noindent \textbf{Corollary 2}: S Suppose \( \mathbf{s}' = (s_1+\frac{n-2}{2}, s_2+\frac{n-2}{2}, \ldots, s_n) \in \mathbb{R}^n \) and  $s_n>-n-1,s_{j}>-\frac{n+3}{2},j = 1,2,\cdots,n-1$. Then we have:

\begin{enumerate}
    \item \( \Delta(y + b) \geq \Delta(y) \), for all \( y, b \in \mathcal{P}_n \);
    \item \( \Delta^{-\mathbf{s}'}(y + b) \leq \Delta^{-\mathbf{s}'}(y) \), for all \( y, b \in \mathcal{P}_n \);
    \item \( |P^{-\mathbf{s}'}(x + iy)| \leq \Delta^{-\mathbf{s}'}(y) \), for all \( x \in \mathbb{R}^m \), \( y \in \mathcal{P}_n \).
\end{enumerate}

\noindent \textbf{Proof}: This corollary can be derived from Lemma 2.3 (2).
 \hfill $\blacksquare$

\begin{lemma} Suppose $\boldsymbol{\eta}=(\eta_1+\frac{n-2}{2},\eta_2+\frac{n-2}{2},\cdots,\eta_n), \boldsymbol{r}=(r_1+\frac{n-2}{2},r_2+\frac{n-2}{2},\cdots,r_n)\in\mathbb{R}^n$ and \(b\in\mathcal{P}_{n}\). If the parameters satisfy the convergence conditions:
$$r_n > \eta_n +\frac{n}{2}+ \frac{1}{2},\quad  \eta_{n}>-1,\quad r_{n}>0, $$
and
$$r_j > \eta_j +n,\quad  \eta_{j}>-\frac{n+1}{2},\quad r_{j}>\frac{1}{2}-\frac{n}{2},$$

then the integral
\begin{align*}
I_{\boldsymbol{r},\boldsymbol{\eta}}(b)
&=\int_{\mathcal{P}_{n}}\Delta^{\boldsymbol{-r}}(y + b)\Delta^{\boldsymbol{\eta}}(y)dy\\
&=C_{5,\mathbf{r},\boldsymbol{\eta}}\Delta^{\eta-\mathbf{r}}(b) \Delta(b)^{\frac{n+1}{2}},
\end{align*}
where
$$C_{5,\mathbf{r},\boldsymbol{\eta}}=C_{4,\mathbf{r}}^{-1} C_{4,\boldsymbol{\eta}} C_{4,\mathbf{r}-\boldsymbol{\eta}}$$

\noindent \textbf{Proof}:
\begin{align*}
I_{\mathbf{r},\boldsymbol{\eta}}(b) &= C_{4,\mathbf{r}}^{-1} \int_{\mathcal{P}_n} \int_{\mathcal{P}_n} e^{-2\pi(y+b)t} \Delta^{-\mathbf{r}}(t^{-1}) \Delta(t)^{-\frac{n+1}{2}} dt \Delta^{\boldsymbol{\eta}}(y) dy \\
&= C_{4,\mathbf{r}}^{-1} \int_{\mathcal{P}_n} e^{-2\pi bt} \Delta^{-\mathbf{r}}(t^{-1}) \Delta(t)^{-\frac{n+1}{2}}dt \int_{\mathcal{P}_n} e^{-2\pi yt} \Delta^{\boldsymbol{\eta}}(y) dy  \\
&= C_{4,\mathbf{r}}^{-1} C_{4,\eta} \int_{\mathcal{P}_n} e^{-2\pi bt} \Delta^{-\mathbf{r}+\boldsymbol{\eta}}(t^{-1}) \Delta(t)^{-n-1} dt \\
&= C_{4,\mathbf{r}}^{-1} C_{4,\boldsymbol{\eta}} C_{4,\mathbf{r}-\boldsymbol{\eta}} \Delta^{\eta-\mathbf{r}}(b) \Delta(b)^{\frac{n+1}{2}}.
\end{align*}
 \hfill $\blacksquare$

\end{lemma}

\begin{lemma}
 Let $\boldsymbol{r}=(r_1+\frac{n-2}{2},r_2+\frac{n-2}{2},\cdots,r_n)\in\mathbb{R}^{n}$, $r_n>\frac{n}{2}+\frac{1}{2}$, $r_j>\frac{3}{2},j = 1,\cdots,n-1$ and $w = u+iv\in T_{\mathcal{P}_n}$.

Then we have
\[
\int_{\mathbb{R}^m} \left|P^{-\mathbf{r}}(u + iv)\right|du = C_{6,\mathbf{r}}\Delta^{-\mathbf{r}}(v)\Delta(v)^{\frac{n + 1}{2}}
\]
Where $$C_{6,\mathbf{r}} = 2^{-(r_1 + r_2+\cdots+r_n)+2n-1}C_{4,\frac{\mathbf{r}}{2}}^{-1}C_{4,\mathbf{r}}.$$

\noindent \textbf{Proof}:
From Tonelli's theorem and Lemma 2.3, we have
\begin{equation}
\begin{aligned}
&\int_{\mathbb{R}^{m}}\vert P^{-\mathbf{r}}(u + iv)\vert du\\
&=C_{4,\mathbf{r}}^{-1}\int_{\mathbb{R}^m}\int_{\mathcal{P}_{n}}e^{2\pi i(u + iv)\cdot t}\Delta^{-\mathbf{r}}(t_{\delta}^{-1})\Delta(t_{\delta})^{-\frac{n + 1}{2}}dtdu\\
&=C_{4,\mathbf{r}}^{-1}\int_{\mathbb{R}^{m}}\int_{\mathbb{R}^m}e^{2\pi i u\cdot t}e^{-2\pi v\cdot t}\Delta^{-\mathbf{r}}(t_{\delta}^{-1})\Delta(t_{\delta})^{-\frac{n + 1}{2}}\chi_{\mathcal{P}_{n}}(t)dtdu.\\
\end{aligned}
\end{equation}

Let
\[f_{v}(t)=e^{-2\pi vt}\Delta^{-\frac{\mathbf{r}}{2}}(t_{\delta}^{-1})\Delta(t_{\delta})^{-\frac{n + 1}{2}}\chi_{\mathcal{P}_{n}}(t).\]

Then
\[f_{v}(t)\in L^{2}(\mathbb{R}^{m}).\]

By Plancherel's theorem, it is proved.
 \hfill $\blacksquare$

\end{lemma}

\begin{lemma}
 If the parameters \(\mathbf{l}\), \(\mathbf{r}\) and \(\boldsymbol{\eta}\) satisfy
$$l_n>-1,r_n>0,\eta_n>\frac{n}{2}+\frac{1}{2},r_n+\eta_n - l_n> n+1,$$ and
$$l_j>-\frac{n+1}{2},r_j>\frac{n-1}{2},\eta_j>n,r_j+\eta_j - l_j>\frac{3n+1}{2}$$
$j = 1,\cdots,n-1.$ Then for \(z,\xi\in T_{\mathcal{P}_{n}}\),
\[
I_{\mathbf{l},\mathbf{r},\boldsymbol{\eta}}(z,\xi)=\int_{T_{\mathcal{P}_n}}\frac{\Delta^{\mathbf{l}}(\text{Im }w)}{P^{\mathbf{r}}(z - \overline{w})P^{\boldsymbol{\eta}}(w - \overline{\xi})}dw = C_{7,\mathbf{l},\mathbf{r},\mathbf{\eta}}P^{-\mathbf{r}-\boldsymbol{\eta}+\mathbf{l}}(z - \overline{\xi})P(z - \overline{\xi})^{n + 1}
\]
where
$$C_{7,\mathbf{l},\mathbf{r},\boldsymbol{\eta}}=C_{4,\mathbf{r}}^{-1}C_{4,\boldsymbol{\eta}}^{-1}
C_{3,\mathbf{l}}C_{4,\mathbf{r}+\mathbf{l}-\boldsymbol{\eta}}.$$

\noindent \textbf{Proof}: Let \(w = u+iv\). By Lemma $2.3$, we have
\begin{equation}
\begin{aligned}
I_{\mathbf{l}, \mathbf{r},\boldsymbol{\eta}}(z,\xi)&=\int_{\mathcal{P}_{n}}\Delta^{\mathbf{l}}(v)\int_{\mathbb{R}^m}P^{-\mathbf{r}}(z - \overline{w})P^{-\boldsymbol{\eta}}(w - \overline{\xi})dudv\\
&=C_{4,\mathbf{r}}^{-1}C_{4,\boldsymbol{\eta}}^{-1}\int_{\mathcal{P}_{n}}\Delta^{\mathbf{l}}(v)\int_{\mathbb{R}^{m}}\left(\int_{\mathcal{P}_{n}}e^{2\pi i(z - \overline{w})\cdot t}\Delta^{-\mathbf{r}}(t_{\delta}^{-1})\Delta(t_{\delta})^{-\frac{n + 1}{2}}dt\right)\\
&\left(\int_{\mathcal{P}_{n}}e^{2\pi i(w-\xi)\cdot t}\Delta^{-\boldsymbol{\eta}}(t_{\delta}^{-1})\Delta(t_{\delta})^{-\frac{n + 1}{2}}dt\right)dudv\\
&=C_{4,\mathbf{r}}^{-1}C_{4,\boldsymbol{\eta}}^{-1}\int_{\mathcal{P}_{n}}\Delta^{l}(v)\int_{\mathbb{R}^{m}}
dudv\\
&\cdot\left(\int_{\mathbb{R}^{m}}e^{-2\pi iut}e^{2\pi izt}e^{-2\pi ivt}\Delta^{-\mathbf{r}}(t_{\delta}^{-1})\Delta(t_{\delta})^{-\frac{n + 1}{2}} \chi_{\mathcal{P}_{n}}(t)dt\right)\\
&\cdot\left(\int_{\mathbb{R}^{m}}e^{2\pi iut}e^{2\pi i\overline{\xi}t}e^{-2\pi ivt}\Delta^{-\mathbf{r}}(t_{\delta}^{-1})\Delta(t_{\delta})^{-\frac{n + 1}{2}} \chi_{\mathcal{P}_{n}}(t)dt\right)
\end{aligned}
\end{equation}
Let
\[f_{z,v}(t)=e^{-2\pi ivt}e^{2\pi izt}\Delta^{-\mathbf{r}}(t_{\delta}^{-1})\Delta(t_{\delta})^{-\frac{n + 1}{2}}\chi_{\mathcal{P}_{n}}(t)\]
\[g_{\xi,v}(t)=e^{-2\pi ivt}e^{2\pi i\overline{\xi} t}\Delta^{-\mathbf{r}}(t_{\delta}^{-1})\Delta(t_{\delta})^{-\frac{n + 1}{2}}\chi_{\mathcal{P}_{n}}(t).\]

Applying Lemma $2.3$, we have \(f_{z, v}(t)\in L^{1}(\mathbb{R}^{m})\) and \(g_{\xi, v}(t)\in L^{1}(\mathbb{R}^{m})\cap L^{2}(\mathbb{R}^{m})\). $(2.2)$ can be expressed as
\[I_{\mathbf{l}, \mathbf{r},\boldsymbol{\eta}}(z,\xi)=C_{4, \mathbf{r}}^{-1}C_{4, \boldsymbol{\eta}}^{-1}\int_{\mathcal{P}_{n}}\int_{\mathbb{R}^{m}}\hat{f}_{z, v}(u)\hat{g}_{\xi, v}(-u)dudv.\]
where \(\hat{f},\hat{g}\) represent the Fourier transforms of \(f, g\) respectively. According to Lemma $2.5$, we have \(\hat{g}_{\xi, v}(-u)\in L^{1}(\mathbb{R}^{n})\). We get
\begin{equation}
\begin{aligned}
I_{\mathbf{l}, \mathbf{r},\boldsymbol{\eta}}(z,\xi)&=C_{2,\mathbf{r}}^{-1}C_{2,\boldsymbol{\eta}}^{-1}\int_{\mathcal{P}_{n}}\Delta^{\mathbf{l}}(v)\int_{\mathbb{R}^m}f_{z,v}(u)\hat{g}_{\xi,v}(-u)dudv\\
&=C_{4,\mathbf{r}}^{-1}C_{4,\boldsymbol{\eta}}^{-1}\int_{\mathcal{P}_{n}}\Delta^{\mathbf{l}}(v)\int_{\mathbb{R}^m}f_{z,v}(u)g_{\xi,v}(u)dudv\\
&=C_{4,\mathbf{r}}^{-1}C_{4,\boldsymbol{\eta}}^{-1}\int_{\mathcal{P}_{n}}\Delta^{\mathbf{l}}(v)\int_{\mathcal{P}_{n}}e^{-2\pi v\cdot u}e^{2\pi i(z - \xi)\cdot u}\Delta^{-\mathbf{r}-\boldsymbol{\eta}}(u_{\delta}^{-1})\Delta(u_{\delta})^{-n - 1}dudv\\
&=C_{4,\mathbf{r}}^{-1}C_{4,\boldsymbol{\eta}}^{-1}\int_{\mathcal{P}_{n}}e^{2\pi i(z - \xi)\cdot u}\Delta^{-\mathbf{r}-\boldsymbol{\eta}}(u_{\delta}^{-1})\Delta(u_{\delta})^{-n - 1}\Delta_{n-1}(u_{\delta})^{n-2}\int_{\mathcal{P}_{n}}e^{-4\pi v\cdot u}\Delta^{\mathbf{l}}(v)dvdu\\
&=C_{4,\mathbf{r}}^{-1}C_{4,\boldsymbol{\eta}}^{-1}C_{3,\mathbf{l}}\int_{\mathcal{P}_{n}}e^{2\pi i(z - \xi)\cdot u}\Delta^{-\mathbf{r}-\boldsymbol{\eta}+\mathbf{l}}(u_{\delta}^{-1})\Delta(u_{\delta})^{-n - 1-\frac{n + 1}{2}}du\\
&=C_{4,\mathbf{r}}^{-1}C_{4,\boldsymbol{\eta}}^{-1}C_{3,\mathbf{l}}C_{4,\mathbf{r}+\mathbf{l}-\boldsymbol{\eta}}P^{-\mathbf{r}-\boldsymbol{\eta}+\mathbf{l}}(z - \overline{\xi})P(z - \overline{\xi})^{n + 1}.
\end{aligned}
\end{equation}
\end{lemma}
 \hfill $\blacksquare$
\begin{lemma}
 If the parameters \(\mathbf{l}\), \(\mathbf{r}\) satisfy
$$l_n>-1, r_n - l_n>n+1,l_j>-\frac{n+1}{2}, r_j - l_j>\frac{3n+1}{2},j = 1,\cdots,n-1.$$

Then
\[
I_{\mathbf{l},\mathbf{r}}(z)=\int_{T_{\mathcal{P}_n}}\frac{\Delta^{\mathbf{l}}(\text{Im }w)}{\left|P^{\mathbf{r}}(z - \overline{w})\right|}dw = C_{8,\mathbf{l},\mathbf{r}}\Delta^{-\mathbf{r}+\mathbf{l}}(y)\Delta(y)^{n + 1}.
\]

Where
\[C_{8,\mathbf{l},\mathbf{r}} = C_{6,\mathbf{r}}C_{5,\mathbf{r}-\mathbf{l}}.\]

\noindent \textbf{Proof}: It can be obtained similarly to Lemma $2.6$.
\end{lemma}
 \hfill $\blacksquare$

\section{Proof of Theorem 1}
 Suppose \(1 < p\leq q<\infty\). Define $\mathbf{R}=(R_1,R_2,\cdots,R_n)$ and $R_j > 0\quad j = 1,2,\cdots,n.$ Let
\[f_{\mathbf{R}}(w)=\frac{\Delta^{\mathbf{l}}(\text{Im}w)}{P^{\mathbf{r}}(w + i\mathbf{R})}, \quad w\in T_{\mathcal{P}_{n}}\]
where $\boldsymbol{l}=(l_1+\frac{n-2}{2},l_2+\frac{n-2}{2},\cdots,l_n)\in\mathbb{R}^{n}$ and $\boldsymbol{r}=(r_1+\frac{n-2}{2},r_2+\frac{n-2}{2},\cdots,r_n)\in\mathbb{R}^{n}$ satisfy the following conditions for all \(j = 1,2,\cdots,n-1\)
\begin{equation}
\begin{aligned}\begin{cases}
&l_{j}>\max\left\{-\frac{n+1}{2p}- \frac{\alpha_{j}}{p},-b_{j} - \frac{1}{2}\right\},\\
&r_{j} > n,\\
&r_{j} -l_{j}>\max\left\{\frac{\alpha_{j}}{p}+\frac{3n+1}{2p},\frac{3n+1}{2}+b_{j} - c_{j}\right\}.\end{cases}
\end{aligned}
\end{equation}

When $j=n$, $l_{n}$ and $r_{n}$ satisfy the following conditions
\begin{align*}\begin{cases}
&l_{n}>\max\left\{-\frac{1}{p}- \frac{\alpha_{n}}{p},-b_{n} - 1\right\},\\
&r_{n} > \frac{n+1}{2},\\
&r_{n} -l_{n}>\max\left\{\frac{\alpha_{n}+n+1}{p},n+1+b_{n} - c_{n}\right\}.
\end{cases}
\end{align*}

We first calculate the norm of \(f_{\mathbf{R}}(w)\) in \(L_{\boldsymbol{\alpha}}^{p}(T_{\mathcal{P}_{n}})\) by $(3.1)$.
\begin{align*}
\left\|f_{\mathbf{R}}(\cdot)\right\|_{L_{\boldsymbol{\alpha}}^{p}}
&=\left(\int_{T_{\mathcal{P}_{n}}}\frac{|\Delta^{p\mathbf{l}}(\mathrm{Im}w)|}{|P^{p\mathbf{r}}(w +i\mathbf{R})|}\Delta^{\alpha}(\mathrm{Im}w)dV(w)\right)^{\frac{1}{p}}\\
&=C_{8,p(\mathbf{l+r})}^{\frac{1}{p}}\Delta^{\mathbf{l}+\frac{\alpha}{p}-\mathbf{r}}(\mathbf{R})
\Delta(\mathbf{R})^{\frac{n+1}{p}}\\
&=C'\prod_{j = 1}^{n}R_{j}^{l_{j}-r_{j}+\frac{\alpha_{j}+n+1}{p}},
\end{align*}
where $C'$ is a constant that depends on the parameters $\boldsymbol{\alpha}$, $\mathbf{l}$, $\mathbf{r}$, $p$ and $q$.

According to Lemma 2.6 and condition (3.1), we have
\begin{align*}
Tf_{\mathbf{R}}(z)&=\Delta^{\mathbf{a}}(\text{Im }z)\int_{T_{\mathcal{P}_{n}}}\frac{\Delta^{\mathbf{b}+\mathbf{l}}(\text{Im }w)}{P^{\mathbf{c}}(z - \overline{w})P^{\mathbf{r}}(w + i\mathbf{R})}dV(w)\\
&=C_{7,\mathbf{b}+\mathbf{l},\mathbf{c},\mathbf{r}}\frac{\Delta^{\mathbf{a}}(\text{Im }z)}{P^{\mathbf{r}+\mathbf{c}-\mathbf{l}-\mathbf{b}}(z + i\mathbf{R})P(z + i\mathbf{R})^{-n - 1}}.
\end{align*}

Since the operator \(T\) is bounded, we have \(\|T_{f_{\mathbf{R}}}\|_{L_{\boldsymbol{\beta}}^{q}}<\infty\). Using Lemma 2.7, we obtain the following conditions: for $j = 1,2,\cdots,n-1$,

\begin{align*}\begin{cases}
qa_j + \beta_j > -\frac{n+1}{2} \\
c_j - b_j - a_j - n - 1 + r_j - b_j-l_j >\frac{\beta_{j}}{q} +\frac{3n+1}{2q},
\end{cases}
\end{align*}

and
\begin{equation}\begin{cases}
qa_n + \beta_n > -1 \\
c_n - b_n - a_n - n - 1 + r_n- b_n-l_n > \frac{\beta_{n}}{q}+\frac{n+1}{q},
\end{cases}
\end{equation}

Moreover,
\[\|Tf_{\mathbf{R}}\|_{L_{\beta}^{q}}=C''\prod_{j = 1}^{n} R_j^{a_j + b_j - c_j + l_j - r_j + n + 1+\frac{\beta_{j}+n+1}{q}}\]
where \(C''\) is a constant depending only on $\mathbf{a}$, $\mathbf{b}$, $\mathbf{c}$, \(\mathbf{l}\), \(\mathbf{r}\), \(\boldsymbol{\beta}\), p and \(q\). Due to the boundedness of the operator \(T\) from \(L_{\boldsymbol{\alpha}}^{p}(T_{\mathcal{P}_{n}})\) to \(L_{\boldsymbol{\beta}}^{q}(T_{\mathcal{P}_{n}})\), we have
\[\|Tf_{\mathbf{R}}\|_{L_{\beta}^{q}}\leqslant C\|f_{\mathbf{R}}\|_{L_{\alpha}^{p}}\]
That is,
\[C''\prod_{j = 1}^{n} R_j^{a_j + b_j - c_j + l_j - r_j + n + 1+\frac{\beta_{j}+n+1}{q}}\leqslant CC'\prod_{j = 1}^{n}R_{j}^{l_{j}-r_{j}+\frac{\alpha_{j}+n+1}{p}}\]
where \(C\), \(C'\) and \(C''\) are independent of \(R_{j}\). For the choice of \(R_{j}\), we only require that it is a positive integer. Therefore, for the above inequality to hold, the following condition must be satisfied:
\begin{equation}
c_j = a_j + b_j + n + 1+\frac{\beta_j+ n + 1}{q}-\frac{\alpha_j + n + 1}{p},\ j = 1, \cdots, n.
\end{equation}
Combining conditions (3.1) and (3.3), condition (3.2) is equivalent to
\begin{equation} -a_jq< \beta_j+\frac{n+1}{2},\ j = 1, \cdots, n-1, \quad -a_nq< \beta_n+1.\end{equation}

We assume that \(\frac{1}{p}+\frac{1}{p'}=1\) and \(\frac{1}{q}+\frac{1}{q'}=1\) . The boundedness of the operator \(T\) from \(L_{\boldsymbol{\alpha}}^{p}(T_{\mathcal{P}_{n}})\) to \(L_{\boldsymbol{\beta}}^{q}(T_{\mathcal{P}_{n}})\) implies the boundedness of its dual operator \(T^{*}\) from \(L_{\boldsymbol{\beta}}^{q'}(T_{\mathcal{P}_{n}})\) to \(L_{\boldsymbol{\alpha}}^{p'}(T_{\mathcal{P}_{n}})\). Note that
\begin{equation}
T^{*}f(z)=\Delta^{\mathbf{b} - \boldsymbol{\alpha}}(\text{Im}z)\int_{T_{\Lambda_{n}}}\frac{\Delta^{\mathbf{a} + \boldsymbol{\beta}}(\text{Im}w)}{P^{\mathbf{c}}(z-\overline{w})}f(w)dV(w).
\end{equation}
Therefore,  we can easily derive that
$$\alpha_j+\frac{n+1}{2} < p\left(b_j+\frac{n+1}{2}\right)\ j = 1, \cdots, n-1,\quad  \alpha_n+1 < p\left(b_n+1\right).$$

\section{Proof of Theorem 2}
We only prove the case when \(j = n\), and the cases for \(j = 1,\ldots,n - 1\) can be proved analogously.
Let
\[H(z,w)=\frac{\Delta^{\mathbf{a}}(\text{Im}z)\Delta^{\mathbf{b} - \alpha}(\text{Im}w)}{\vert P^{\mathbf{c}}(z-\overline{w})\vert},\]
\[\varphi_{1}(w)=\Delta^{\mathbf{r}}\text{Im}w,\quad\varphi_{2}(z)=\Delta^{\mathbf{l}}(\text{Im}z).\]

Here, the indices \(\mathbf{r}\) and \(\mathbf{l}\) need to satisfy certain conditions so that we can use Lemma 2.2. By Lemma 2.2, let \(\frac{1}{p}+\frac{1}{p'}=1\)

\begin{equation}
\begin{aligned}
&\int_{T_{\mathcal{P}_{n}}} H(z, w)^{tp'} \phi_1(w)^{p'} dV_{\alpha}(w)\\
=&\int_{T_{\mathcal{P}_{n}}} \left(\frac{\Delta^{\mathbf{a}}(\text{Im }z)\Delta^{\mathbf{b}-\mathbf{\alpha}}(\text{Im }w)}{\vert P^{\mathbf{c}}(z - \overline{w})\vert}\right)^{tp'} \Delta^{p'r + \alpha}(\text{Im }w)dV(w)\\
=&\Delta^{tp'\mathbf{a}}(\text{Im }z)\int_{T_{\mathcal{P}_{n}}} \frac{\Delta^{tp'(\mathbf{b}-\mathbf{\alpha})+p'r + \alpha}(\text{Im }w)}{\vert P^{tp'\mathbf{c}}(z - \overline{w})\vert}dV(w)\\
=&C_1\Delta^{tp'\mathbf{b}-tp'\mathbf{\alpha}+p'r + \alpha - tp'\mathbf{c}+tp'\mathbf{a}}(\text{Im }z)\Delta(\text{Im }z)^{n + 1}\\
\end{aligned}
\end{equation}
where \( C_1 \) is a positive constant. Similarly, for
\begin{equation}
\begin{aligned}
&\int_{T_{\mathcal{P}_{n}}}H(z,w)^{q(1 - t)}\varphi_{1}(z)^{s}dV_{\beta}(z)\\
=&\int_{T_{\mathcal{P}_{n}}}\left(\frac{\Delta^{\mathbf{a}}(\text{Im}z)\Delta^{\mathbf{b} - a}(\text{Im}w)}{\vert P^{\mathbf{c}}(z-\overline{w})\vert}\right)^{q(1 - t)}\Delta^{ql+\beta}(\text{Im}z)dV(z)\\
=&\Delta^{q(1 - t)(b - \alpha)}(\text{Im}w)\int_{T_{\mathcal{P}_{n}}}\frac{\Delta^{q(1 - t)a+ql+\beta}(\text{Im}z)}{\vert P^{q(1 - t)c}(z-\overline{w})\vert}dV(z)\\
=&C_2\Delta^{q(\mathbf{a}-\mathbf{c}+\mathbf{b}-\mathbf{\alpha})+tq(\mathbf{c}-\mathbf{a}+\mathbf{b}-\mathbf{\alpha})+q\mathbf{l}+\beta}(\text{Im }w)\Delta(\text{Im }w)^{n + 1}\\
\end{aligned}
\end{equation}
where \( C_2 \) is a positive constant. Let \( c' = \min\{c_1, c_2, \cdots, c_n\} \). We choose \( t \) such that
\begin{equation}
\frac{n}{c'} - \frac{n}{p c'} < t < 1 - \frac{n}{q c'}.
\end{equation}
Since for each \( j \), \( c_j > n \), we have
\[
1 - \frac{n}{q c'} - \frac{n}{c'} + \frac{n}{p c'} > \frac{n}{p c'} - \frac{n}{q c'} \geq 0
\]
and \( 0 < \frac{n}{c'} - \frac{n}{p c'} < 1 \) and \( 0 < 1 - \frac{n}{q c'} < 1 \). This ensures the existence of such \( t \) and the chosen \( t \) must satisfy \( 0 < t < 1 \) and
\[
\max_j \left\{ \frac{n}{c_j} - \frac{n}{p c_j} \right\} < t < \min_j \left\{ 1 - \frac{n}{q c_j} \right\}.
\]
Next, we let
\begin{align*}\begin{cases}
A_n = - \alpha_n + \frac{\alpha_n}{p} + t \alpha_n- t b_n\\
B_n = t(c_n - b_n + \alpha_n) - \alpha_n + \frac{\alpha_n}{p} - n - 1 + \frac{n+1}{p}\\
C_n = \frac{n }{q} + (t - 1)(c_n - b_n + \alpha_n)\\
D_n =- (1 - t)\alpha_n + (1 - t)b_n.
\end{cases}
\end{align*}

(4.3) ensures that \( A_n < B_n \) and \( C_n < D_n \), i.e., the intervals \( (A_n, B_n) \) and \( (C_n, D_n) \) are non-empty. $\alpha_n + 1 < p( b_n + 1)$ and $c_j = a_j + b_j + n + 1 + \frac{\beta_j + n + 1}{q} - \frac{\alpha_j + n + 1}{p}
$ ensure \( C_n < B_n \), while $\beta_n + 1 > -qa_n$ and $c_j = a_j + b_j + n + 1 + \frac{\beta_j + n + 1}{q} - \frac{\alpha_j + n + 1}{p}
$ imply \( A_n < D_n \). Therefore, we have

\[
(A_j, B_j) \bigcap (C_j, D_j) \neq \varnothing.
\]
Thus, we choose \( r_n \in (A_n, B_n) \bigcap (C_n, D_n) \neq \varnothing \). Let \( 1/p + 1/p' = 1 \), then
\begin{equation}
\begin{cases}
t p'(b_n - \alpha_n) + p' r_n + \alpha_n > -1 \\
t p' c_n - t p'(b_n - \alpha_n) - p' r_n - \alpha_n > n + 1.
\end{cases}
\end{equation}
Condition (4.4) is derived from \( r \in (A_n, B_n) \). We set
\[
l_n = r_n + c_n - b_n - a_n + \alpha_n - t(c_n - b_n - a_n + \alpha_n) - \frac{\beta_n+ n + 1}{q}.
\]
Then, from \( r_n \in (C_n, D_n) \), we can obtain that \( l_n \) satisfies
\begin{equation}
\begin{cases}
q(1 - t)a_n + q l_n + \beta_n > -1 \\
q(1 - t)c_n - (q(1 - t)a_n + q l_n + \beta_n) > n+1.
\end{cases}
\end{equation}

From (4.1), $c_j = a_j + b_j + n + 1 + \frac{\beta_j + n + 1}{q} - \frac{\alpha_j + n + 1}{p}$ and the definition of \( l_n \), we have
\[
-t p'(c_n - b_n - a_n + \alpha_n) + p' r_n + \alpha_n + n + 1 = p' l_n
\]
Hence
\[
\int_{T_{\mathcal{P}_n}} H(z, w)^{t p'} \phi_1(w)^{p'} dV_\alpha(w) = C_1 \Delta^{p'l}(\text{Im} z) = C_1 \phi_2(z)^{p'}.
\]

From (4.2) and $c_j = a_j + b_j + n + 1 + \frac{\beta_j + n + 1}{q} - \frac{\alpha_j + n + 1}{p}$, we have
\[
q(a_n - c_n + b_n - \alpha_n) + t q(c_n - a_n - b_n + \alpha_n) + q l_n + \beta_n+ n + 1 = q r_n.
\]

Thus, we have
\[
\int_{T_{\mathcal{P}_n}} H(z, w)^{q(1-t)} \phi_1(z)^q dV_\beta(z) = C_2 \Delta^{q r}(\text{Im} w) = C_2 \phi_1(w)^q.
\]
By Lemma 2.2, we can conclude that the operator \( T \) is bounded from \( L_\alpha^p(T_{\mathcal{P}_n}) \) to \( L_\beta^q(T_{\mathcal{P}_n}) \).

\hfill $\blacksquare$


\begin{thebibliography}{999}

\bibitem{1} Bergman S. The kernel function and conformal mapping[M]. Providence: American Mathematical Society, 1950.
\bibitem{2} Deng G T, Huang Y, Qian T. Reproducing kernels of some weighted Bergman spaces[J]. The Journal of Geometric Analysis, 2021, 31(10): 9527 - 9550.
\bibitem{3} Deng G T, Fu Q, Cao H. Laplace transforms for analytic functions in tubular domains[J]. Acta Mathematica Scientia, 2021, 41(6): 1938 - 1948.
\bibitem{4} Forelli F, Rudin W. Projections on spaces of holomorphic functions in balls[J]. Indiana University Mathematics Journal, 1974, 24(6): 593 - 602.
\bibitem{5} Hedenmalm H, Korenblum B, Zhu K H. Theory of Bergman spaces[M]. Springer Science Business Media, 2012.
\bibitem{6} Kures O, Zhu K H. A class of integral operators on the unit ball of $\mathbb{C}^n$[J]. Integral Equations and Operator Theory, 2006, 56(1): 71 - 82.
\bibitem{7} Kolaski C J. A new look at a theorem of Forelli and Rudin, Indiana Univ. Math. J. 28(1979),495-499.
\bibitem{8} Li S X, Wulan H, Zhu K H. A characterization of Bergman spaces on the unit ball of $\mathbb{C}^{n}$[J]. Canadian Mathematical Bulletin- Bulletin Canadien de Mathematiques, 2012, 55(1): 146-152.
\bibitem{9} Liu C W, Liu Y, Hu P Y, et al. Two classes of integral operators over the Siegel upper half-space[J]. Complex Analysis and Operator Theory, 2019, 13(3): 685 - 701.
\bibitem{10} Liu C W, Si J J. Positive Toeplitz operators on the Bergman spaces of the Siegel upper half-space[J]. Communications in Mathematics and Statistics, 2020, 8(1): 113 - 134.
\bibitem{11}Liu J X, Deng G T, Gao Z Q. Boundedness of Forelli-Rudin Type Operators on Tube Domains over the Forward Light Cones [J]. Acta Mathematica Scientia, 2025, 45B(2): 1-13.
\bibitem{12} McNeal. The Bergman projection as a singular integral operator[J]. The Journal of Geometric Analysis, 1994, 4: 91 - 103.
\bibitem{13}Nana C, Sehba B F. Off-diagonal estimates of some Bergman-type operators of tube domains over symmetric
cones. Positivity, 2018, 22: 507-531
\bibitem{14}Sehba B F. Bergman type operators in tubular domains over symmetric cones. Proceedings of the Edinburgh Mathematical Society, 2009, 52(2): 529-544
\bibitem{15}Stein E M. Singular integrals and estimates for the Cauchy-Riemann equations, Bull. Amer. Math. Soc. 79(1973), 440-445.
\bibitem{16} Wang  X, Liu M S. The boundedness of two classes of integral operators[J]. Czechoslovak Mathematical Journal, 71(2), 475 - 490.
\bibitem{17} Yang K. Essential norms of weighted composition operators on weighted Bergman spaces over the unit ball[J]. Applied Mathematics (English Edition), 2010, 2010(1): 41 - 48.
\bibitem{18} Zhu K H. A Forelli - Rudin type theorem with applications[J]. Complex Variables, Theory and Application: An International Journal, 1991, 16(2-3): 107 - 113.
\bibitem{19} Zhu K H. A sharp norm estimate of the Bergman projection on $L^{p}$ spaces[J]. Contemporary Mathematics, 2006, 404: 199.
\bibitem{20} Zhu K H. Spaces of holomorphic functions in the unit ball[M]. New York: Springer, 2005.
\bibitem{21} Zhu K H. Embedding and compact embedding of weighted Bergman spaces[J]. Illinois Journal of Mathematics, 2022, 66(3): 435 - 448.
\bibitem{22}Zhao R H. Generalization of Schur's test and its application to a class of integral operators on the unit ball of $C^{n}$[J]. Integral equations and operator theory, 2015, 82(4): 519-532.
\bibitem{23} Zhao R H, Zhou L F. $L^p - L^q$ boundedness of Forelli-Rudin type operators on the Unit ball of $\mathbb{C}^n$[J]. Journal of Functional Analysis, 2022, 282(5): 109345.
\bibitem{24} Zhao R H, Zhu K H. Theory of Bergman spaces in the unit ball of $\mathbb{C}^n$[J]. Memoires de la Societe mathematique de France, 2006, 115.
\bibitem{25} Zhou L F, Wang X, Liu M S. The boundedness of Forelli-Rudin type operators on the Siegel upper half-space[J]. Complex Analysis and Operator Theory, 2023, 17(8): 127.


\end{thebibliography}
\end{document}